

\input amstex
\documentstyle{amsppt}

\magnification\magstep1

\define\p{\Bbb P}

\define\a{\Bbb A}

\redefine\c{\Bbb C}

\redefine\o{\Cal O}

\define\q{\Bbb Q}

\define\r{\Bbb R}

\define\z{\Bbb Z}

\define\n{\Bbb N}

\define\map{\dasharrow}

\define\qtq#1{\quad\text{#1}\quad}

\define\section#1{

\bigpagebreak{\smc#1}\bigpagebreak

}

\define\rup#1{\ulcorner #1\urcorner}

\define\rdown#1{\llcorner #1\lrcorner}

\define\demop{\demo{Proof}}

\define\pic{\operatorname{Pic}}

\define\deq{:=}

\define\norm#1{\vert\vert#1\vert\vert}

\define\broot#1{[@,@,@,\root{n}\of{#1}@,@,@,]}

\define\const{\text{const}}

\rightheadtext{Low degree polynomial equations}

\topmatter
\title     Low degree polynomial equations: arithmetic, geometry and topology
\endtitle
\author
J\'anos Koll\'ar
\endauthor
 \endtopmatter

Polynomials appear in mathematics frequently, and we all know from
experience that   low degree polynomials are easier to deal with than high
degree ones.  It is, however, not clear that there is a well defined class of
``low degree" polynomials. For many questions,
polynomials behave well if their degree is low enough, but the
precise bound on the degree depends on the concrete problem.

My interest is to investigate polynomials through their zero sets. That is,
using sets of the form
$$
\{(x_1,\dots,x_n)\vert f(x_1,\dots,x_n)=0\}.
$$
  I   intentionally refrain from specifying where the
coordinates $x_i$ are. They could be rational, real or complex numbers, but in
some cases the $x_i$ will be polynomials in a new variable $t$. My 
focus is on the polynomial $f$.

Consider, for instance, a polynomial 
$$
f:= a_0+\sum_{i=1}^na_ix_i^k,\qtq{where}a_i\in \z\setminus\{0\}.
$$
Specifying where the coordinates are, leads us to various branches of
mathematics:

\demo{Arithmetic} Choose $x_i\in \q$. The solutions of these Fermat-type
equations have been  much studied, some cases going back to Diophantus, but we
still know very little if $n>2$.
\enddemo

\demo{Topology} Choose $x_i\in \r$ or $x_i\in \c$. The set of solutions is a
topological manifold, and various topological properties  can be related to
algebraic properties of $f$.  For instance, the dimension  and the homology
can  be computed in terms of $n,k$. (Over $\r$ we also need to know the signs
of the $a_i$.)
\enddemo

\demo{Complex manifolds} Choose  $x_i\in \c$. The set of
solutions is a complex analytic manifold. The holomorphic function theory of
this complex manifold can  be understood in terms of polynomials. This
is especially true in the compact versions of this problem.  
\enddemo

\demo{Finite fields} We can also look at solutions of $f=0$ in finite fields.
Centuries ago this was   done by studying $f\equiv 0 \mod p$.
Recently, algebraic geometry over finite fields found many connections with  
coding theory, combinatorics and    computer science.
\enddemo

I like to think of any of the zero sets as a snapshot of the  polynomial $f$.
They all show something about $f$. Certain snapshots reveal more than
others:

\demo\nofrills{Do zero sets determine a polynomial?\usualspace}
For instance, $x_1^{2k}+\dots+x_n^{2k}+1=0$  has no solutions in $\q$, not
even in $\r$.  Thus the zero set gives essentially no information.
The situation is very different over algebraically closed
fields. 
If $f,g\in \c[x_1,\dots,x_n]$,
then
$$
\{\bold x\in \c^n\vert f(\bold x)=0\}=
\{\bold x\in \c^n\vert g(\bold x)=0\}
\qquad\Leftrightarrow\qquad
\matrix
\text{$f$ and $g$ have the same}\\
\text{irreducible factors.}
\endmatrix
$$
(This is an easy special   case
of the   Nullstellensatz of \cite{Hilbert1893}.) If we want to go
further, we must study  solutions of $f=0$ in  any  commutative ring   $R$
with a unit. This  approach was first adopted by Grothendieck in
\cite{EGA60-67}, though in retrospect,  \cite{Weil46} and
\cite{Rilke30,vol.2.p.175} clearly pointed in this direction.
 We obtain that if $f,g\in \z[x_1,\dots,x_n]$ are two polynomials, then
$$
\matrix
\{\bold x\in R^n\vert f(\bold x)=0\} =
\{\bold x\in R^n\vert g(\bold x)=0\}\\
\text{(for every commutative ring $R$)}
\endmatrix
\qquad\Leftrightarrow\qquad
f=\pm g.
$$
  Thus studying solutions in all  commutative rings  determines the polynomial
up to a sign. This approach is very powerful, but rather technical. Therefore
I will stick to studying solutions in fields for the rest of the lecture.
\enddemo

It   turns out that there is a collection of basic questions in arithmetic,
algebraic geometry and topology all of which give the same class of ``low
degree" polynomials.  The aim of this lecture is to explain these properties
and to provide a survey of the known results.

\head 1. Introductory Remarks
\endhead

We start with the observation that in some cases the degree alone does not
provide a good measure of the complexity of a polynomial equation. In order to
develop the correct picture, we have to understand which polynomials behave in
an atypical manner.

\demo{1.1 High degree polynomials that behave like low degree ones} 

There are at least three situations when the zero set of a high degree
polynomial shares some of the properties of zero sets of low degree
polynomials:

\demo{1.1.1 Reducible equations}   If $f=gh$, then the set $(f=0)$ is the union
of the sets $(g=0)$ and $(h=0)$. Thus we can   restrict ourselves to
 the case when $f$ is irreducible. 
\enddemo

\demo{1.1.2 Low degree in certain variables} 
Let us consider an extreme case, when $f$ has degree 1 in   the
variable $x_n$. Then $f$ can be written as 
$$
f=f_1(x_1,\dots,x_{n-1})+x_nf_2(x_1,\dots,x_{n-1}).
$$
The substitution $x_n=-f_1/f_2$ shows that the set $(f=0)$ behaves  like the  
vector space of the first $(n-1)$ variables
$\{(x_1,\dots,x_{n-1})\}$.  This is completely true if $f$ is linear, but in
general 
the correspondence  breaks down if 
$f_2=0$. The latter equation involves one fewer variable, and therefore
it is considered easier.  Roughly speaking,  $f$ should be viewed as
complicated as a linear equation. In general, if $f$ has low degree in certain
variables then it behaves like a low degree equation.
\enddemo

\demo{1.1.3  Very singular equations}  Consider for instance the equation
$x_1^d-x_2^{d-1}=0$. Its degree in both variables is high. Nonetheless, the
substitution
$$
x_1=t^{d-1},\ x_2=t^d
$$
shows that solutions of $x_1^d-x_2^{d-1}=0$ are parametrized by  the values of
the variable $t$.  The same happens for any polynomial $f(x_1,x_2)$  of degree
$d$ all of whose partials up to order $d-2$ vanish at a certain point.
In general,  a high degree equation  $f$  
behaves as a low degree equation if many of   the partial derivatives of $f$  
simultaneously vanish  at  many points. 
\enddemo

While all of these cases do occur, there are relatively few polynomials that
behave this way.  For
instance, all  polynomials $f(x_1,x_2)$ of degree $\leq d$ form a
vectorspace
$V_d$ of dimension $\binom{d+2}{2}$. The set of polynomials which are
exceptional for any of the above 3 reasons is  a subset of codimension
$d-1$ for $d\geq 2$. 
\enddemo

 This remark shows that for most polynomials
the degree is a good measurement of complexity. 
In order to run computer experiments, it is desirable to have a class of
polynomials with very few nonzero coefficients which are nonetheless
``general". A good set of examples to keep in mind is the following.

\demo{1.2 Test Examples}  The equations
$\sum_i c_ix_i^{d} =c_0$ have been much  studied. Unfortunately, they are
 sometimes too special. 
It seems that the inhomogeneous version is much more indicative of the general
case.
 Fix natural numbers  $d_i:\ i=1,\dots,n$ and $c_0,\dots,c_n$ such that
$\prod_ic_i\neq 0$. Then
$$
\sum_{i=1}^n c_ix_i^{d_i} =c_0 \qtq{has ``low degree" iff} \sum_{i=1}^n
\frac1{d_i}\geq 1.
\tag 1.2.1
$$
We see in (5.5) that the above condition does coorespond to the eventual
definition (4.1). Moreover, I claim that the behaviour of these examples
correctly predicts the broad features of the theory. You have to trust me that
this  purely experimental assertion is valid.
\enddemo

As a first example, let us see what a simple minded constant count gives  
about solutions of the   equations (1.2.1) over $\q$.

\demo{1.3 Heuristic claim}  Fix natural numbers  $d_i\:i=1,\dots n$ and
rational numbers $c_i\:i=0,\dots n$.  I claim that usually 
$$
\sum_{i=1}^n c_ix_i^{d_i} =c_0 \qtq{ has many solutions in $\q$ iff}
\sum_{i=1}^n
\frac1{d_i}\geq 1.\tag 1.3.1
$$
Unfortunately there are large classes of equations where this is false.  For
instance, $\sum x_i^2=-1$  has
no solutions in $\q$, not even in $\r$. Looking at $x_1^2 -x_2^2$   modulo 
$4$, we see that $x_1^2 -x_2^2=2$ has no rational solutions. There are several
approaches to correct these problems; we   encounter two of them
later.
For the moment I ignore these counterexamples, and give a proof of
(1.3.1).

It is easier to look for
integral solutions, so we homogenize the equation in  the following
(somewhat unusual) way. Set $d_0$
to be the least common multiple of
$d_1,\dots,d_n$ and  let $d_0=d_ib_i$. Look at the equation
$$
\sum_{i=1}^n c_iy_i^{d_i} =c_0y_0^{d_0}.\tag 1.3.2
$$
There is a correspondence between solutions of (1.2.1) and of (1.3.2) given by
$$
(x_1,\dots,x_n)\mapsto (1,x_1,\dots,x_n)\qtq{and}
(y_0,y_1,\dots,y_n)\mapsto  (y_1/y_0^{b_1},\dots,y_n/y_0^{b_n}).
$$
This shows that finding all rational solutions of (1.2.1) is equivalent to
finding all integral solutions of (1.3.2).

Set $f=-c_0y_0^{d_0}+\sum_{i=1}^n c_iy_i^{d_i}$. There is a constant  $C$,
depending on $f$, such that
$$
 |f(y_0,\dots,y_n)|\leq  C\cdot (\max_i|y_i|^{d_i}).\tag 1.3.3
$$
Fix  a large   $N$ and let the $y_i$ run through the set of integers in 
$[-N^{1/d_i},N^{1/d_i}]$. We get 
$$
\const\cdot N^{\sum_{i=0}^n (1/d_i)}
\qtq{values of $f$   in the interval} [-C\cdot N, C\cdot N].
$$ 
{\bf If} these values
are uniformly distributed,   we   obtain the asymptotic 
$$
\#\{\sum_i c_iy_i^{d_i} =c_0y_0^{d_0}, |y_i|\leq N^{1/d_i}\}\sim 
\const\cdot N^{-1+\sum_{i=0}^n (1/d_i)}\qtq{as $N\to\infty$.}
$$
If
$\sum_{i=1}^n (1/d_i)\geq 1$ then 
$\sum_{i=0}^n (1/d_i)> 1$, and
the number of solutions grows as a power of
$N$. 
If $\sum_{i=1}^n (1/d_i)<1$ then $\sum_{i=0}^n (1/d_i)\leq 1$ because of
the special choice of $d_0$, thus there should be  few
solutions. \qed

For which other polynomials $f$ does this counting method work? The main part
is the estimate (1.3.3). This works if $f$ is weighted homogeneous  of degree 1
with weights
$1/d_i$. That is, if  we declare $\deg x_i=1/d_i$ then $\deg f\leq 1$.

There are some examples where the above simpleminded counting method does
work, for instance,   for equations of the form
$$
f(x_1,\dots,x_n)-f(y_1,\dots,y_n)=0.
$$
The above argument gives a lower bound
$$
\#\{f(x_1,\dots,x_n)=f(y_1,\dots,y_n), |x_i|, |y_i|\leq N\}\geq \const\cdot
N^{2n-d}.
$$
This is  interesting only if $d<n$ since   the trivial
solutions $x_i=y_i$ always give a lower bound $\const\cdot N^n$.
\enddemo

In the rest of the lecture I aim to explain the various properties that lead
to this class of equations, starting with the 2-variable case in section 2. 
This is called the theory of
algebraic curves. Most of the theory was well-established in the 19th
century, with the exception of the arithmetic aspects.

Section 3 is devoted to the 3-variable case, which corresponds to the theory of
algebraic surfaces. The geometric aspects have been established around the
turn of the century, many of the topological results are recent and most of
the arithmetical questions are open. 

Much less is known in higher dimensions. The open questions involve deep
problems in algebraic geometry, number theory and differential  topology. I
am confident that these problems constitute a very interesting direction of
research for a long time to come.

\head 2. Two Variable Polynomials = Algebraic Curves
\endhead

Let us consider a 2 variable polynomial $f(x,y)=\sum a_{ij}x^iy^j$ of degree
$d$.  Let $C_{\text{\it aff}}$ denote its zeros, that is,
$$
C_{\text{\it aff}}\deq \{(x,y)\vert f(x,y)=0\}.
$$
(The subscript {\it aff} refers to the fact that we are in affine 2-space
$\a^2$.) This is not a  set since I have not specified where the coordinates
$x,y$ are. If the coefficients $a_{ij}$ are in a field $F$, then for
any larger field $E\supset F$ we can look at  solutions of $f=0$ in $E$. The
resulting set is 
$$
C_{\text{\it aff}}(E)\deq \{(x,y)\in E^2\vert f(x,y)=0\}\subset E^2.
$$
A common case is when $a_{ij}\in \q$, 
and for the larger
field $E$ we choose $\q,\r$ or $\c$.

$C_{\text{\it aff}}(\q)$ is just a set of points, but $C_{\text{\it
aff}}(\r)\subset
\r^2$ naturally appears as a curve (that is, a 1-dimensional topological
space).
$C_{\text{\it aff}}(\c)\subset \c^2$ is a Riemann surface:  a complex
manifold  locally like
$\c$.

In studying the manifolds $C_{\text{\it aff}}(\r)$ or $C_{\text{\it aff}}(\c)$
it is frequently inconvenient that they are not compact. The usual way to deal
with this problem is to introduce the projective plane $\p^2$ with homogeneous
coordinates
$(x_0:x_1:x_2)$. Its relationship to the old affine coordinates is
$x=x_1/x_0, y=x_2/x_0$. 
If the coordinates $x_i$ are in a field   $E$, we obtain the
corresponding projective plane $\text{E}\p^2$. The most frequently used ones
are
$\q\p^2,\r\p^2$ and $\c\p^2$.

The homogenization of $f$ is given by
$$
\bar f(x_0,x_1,x_2)\deq x_0^df(x_1/x_0,x_2/x_0).
$$
The corresponding zero set
$$
C(E)\deq \{(x_0:x_1:x_2)\in \text{E}\p^2\vert \bar f(x_0,x_1,x_2)=0\}\subset
\text{E}\p^2
$$
turns out to be more convenient for most purposes.

Based on the real picture, algebraic geometers say that $C$ is an
{\it algebraic curve}.  Thus we prefer to call $\c$ the
complex line (the complex plane is of course $\c^2$).
This leads to occasional confusion, but this is not the time to change 150
year-old terminology.

In what follows I   collect certain properties of algebraic curves defined by
equations of degree at most 2. In all cases I would like the properties to
hold only for curves defined by equations of degree at most two (assuming the
genericity conditions of (1.1)). 

All of the characterizations listed here are standard results of the theory
of algebraic curves and Riemann surfaces.  
One of the most accessible   introductions to algebraic geometry  is
\cite{Shafarevich94} (or any of the other editions). For algebraic curves see
\cite{Fulton69}. The analytic theory of Riemann surfaces is treated in
\cite{Siegel69; Gunning76}. For the arithmetic aspects I found
\cite{Serre73; Silverman86} especially useful.

\subhead  Characterizations of ``low degree " curves
\endsubhead

I start with the algebraic geometry condition, not because it is the most
obvious for curves, but because this provides the neatest definition in higher
dimensions.

\proclaim{2.1  Algebraic geometry}  There is a one-to-one map
 given by rational functions
$g:\c\p^1\to C(\c)$.
\endproclaim

In this case   $C$ is called a {\it rational curve}. 

Let $(s:t)$ be the homogeneous coordinates on $\c\p^1$.
If $f= a_0x_0+a_1x_1+a_2x_2$ is linear and $a_2\neq 0$,   we can choose
$$
g:(s:t)\mapsto \left(a_2s:a_2t:-(a_0s+a_1t)\right).
$$
For $\deg f=2$  assume for simplicity that 
$f=a_0x_0^2+a_1x_1^2+a_2x_2^2$. (This can always be achieved after a
linear change of coordinates.) We can take
$$
g:(s:t)\mapsto
\left(a_1s^2-a_0t^2: -2a_0st: \sqrt{-a_0/a_2}(a_1s^2+a_0t^2)
\right).
$$
(In case you wonder where this came from, let 
$h:C\to L$ be the projection of $C$ from the point 
$P=(\sqrt{a_2}:0:\sqrt{-a_0})\in C$ to the $(x_2=0)$ line
(Mercator projection). $g$ is the inverse of $h$.)

The fact that no such $g$ exists for higher degree equations is harder.

\proclaim{2.2  Topology}  $C(\c)$ is homeomorphic to the sphere $S^2$.
\endproclaim 

The maps $g$ from (2.1) also provide a homeomorphism; the hard part is again
to see   that this cannot be done for higher degree equations.
The precise result is that if $C$ is defined by a degree $d$ equation then
$C(\c)$ is homeomorphic to a sphere with $\frac12(d-1)(d-2)$ handles.

\proclaim{2.3  Hard Arithmetic}  $C(\q)$ is ``large"
\endproclaim

For this to make sense, we should start with a curve 
$$
C=(\bar f(x_0,x_1,x_2)=0)\subset \p^2,
$$
 where $\bar f$ has rational coefficients. 

Unfortunately it is not easy to pin down what ``large" exactly means. First of
all, if $n\geq 4$ then $C(\q)$ is finite  by \cite{Faltings83}. Unfortunately,
$C(\q)$ is often infinite for $n=3$ and frequently empty for $n=2$.

To get the right answer, we have to develop a  good measure of the size of a
solution. This is most conveniently done in projective coordinates.

Any point   $P\in \q\p^2$ can be represented as a triple $P=(x_0:x_1:x_2)$
where
$x_0,x_1,x_2\in \z$ are relatively prime. This representation is unique up to
sign, thus $H(P)\deq \max\{|x_0|,|x_1|,|x_2|\}$ is well-defined. It is called
the {\it height} of $P$.
One    defines the counting function
$$
N(C,H)\deq \#\{P=(x_0:x_1:x_2)\in \q\p^2\vert \bar f(x_0,x_1,x_2)=0\text{ and }
H(P)\leq H\}.
$$
Roughly speaking, we look for rational solutions of $f(x,y)=0$ where the
numerators and denominators are bounded.

This nearly gives the right answer. If $n=2$ then  $C(\q)$ is either
empty or $N(C,H)$ grows like
$\const\cdot H$; if $n=3$ then $N(C,H)$ grows slower than any power of $H$
\cite{N\'eron65}.

In order to deal with the case when  $C(\q)$  is
empty, we have to count
solutions in various algebraic number fields.
It is not hard to generalize the notion of height to the case when the
coordinates of $P$ are in an algebraic number field $E\supset \q$
(see  \cite{Silverman86,VIII.5} for a short and clear summary). 
 We obtain a
similar counting function $N_E(C,H)$. This finally gives the correct
generalization:

\demo{2.3.1 Theorem} $C$ is a rational curve iff  $N_E(C,H)$ grows polynomially
with
$H$ for a suitable number field $E$.
\enddemo

\proclaim{2.4  Complex manifolds}  $C(\c)$ has genus zero.
\endproclaim

Global holomorphic differential forms on a compact Riemann surface  have been
much studied, starting with the works of Euler, Abel and Riemann. On a Riemann 
surface we   have only
$1$-forms, these are locally given as 
$f(z)dz$
where  $z$ is a local coordinate and  $f(z)$ is
holomorphic.   Such forms are automatically closed, thus the integral
$$
\int_{\gamma} f(z)dz \qtq{over a closed loop}\gamma\subset C(\c)
$$
depends only on the homology class  $[\gamma]\in H_1(C(\c),\z)$. Since the
fundamental studies of Riemann, these  give the basic approach to  finer
understanding of Riemann surfaces.

By definition, the {\it genus} is the dimension of the vector space of 
global holomorphic differential forms. If there are no such forms, the above
integrals give no information. Fortunately, this happens precisely when other
descriptions are very simple.

\proclaim{2.5  Easy Arithmetic} There are many solutions over function fields.
\endproclaim

Here we look at the behaviour of the sets $C(F)$
where $F=\c(t)$ is the field of rational functions in one variable. 
Of course $f=\sum a_{ij}(t)x^iy^j$ and the coefficients $a_{ij}(t)$
themselves are rational functions. The field
$\c(t)$ shares many properties of $\q$, but the results are  easier to state
and the proofs are much simpler. (The difference between $\q$ and
$\c(t)$ becomes apparent when studying their Galois cohomology.)

The advantage of $\c(t)$ is that there are two ways of looking at solutions
over $\c(t)$.

(2.5.1.1) The algebraic way. Just handle everything as quotients of polynomials
in $\c[t]$.

(2.5.1.2) The geometric way. An equation $f(x,y)=0$ with coefficients in
$\c(t)$ can be viewed as an equation $\tilde f(x,y,t)=0$ with coefficients in
$\c$. This defines an algebraic surface $S\subset \c^3$ and we have a
distinguished coordinate projection to the $t$-axis $p:S\to \c_t$. 

A solution $(x(t),y(t))$ of $f(x,y)=0$ can be identified with a map
$$
h:\c_t\to S\qtq{given by} t\mapsto (x(t),y(t),t).
$$
$h$ is a section of $p:S\to \c_t$ and every (rational) section arises as
above.

The first indication that we can expect nicer results is the following
theorem, which can be proved by a straightforward generalization of the
counting argument (1.3).
The first proof is in \cite{Noether1871}. Later algebraic proofs, more
suited to generalizations, are in \cite{Baker22,vol.VI.p.147} and
\cite{Tsen36}.

\demo{2.5.2 Theorem}  If $\deg f\leq 2$ then $f=0$ has a solution in $\c(t)$.
\enddemo

We may also want  to know that there are many solutions. A natural approach
is to look for solutions $(x(t),y(t))$ where certain values
$(x(t_k),y(t_k))$ are specified in advance. This is possible only if
the points $(x(t_k),y(t_k),t_k)$ lie on the surface $S$, that is, if
$\sum a_{ij}(t_k)x(t_k)^iy(t_k)^j=0$. In this case we say that
the pair $(x(t_k),y(t_k))$ is a solution of $f(x,y)=0$ at $t_k$.

As an easy exercise in the theory of algebraic surfaces
we get a very strong
characterization:

\demo{2.5.3 Theorem}  There is a finite set
$B\subset \c$ such that if $t_1,\dots,t_s\in \c\setminus B$ are arbitrary
points and  $(x_k,y_k)$ any solution of $f$ at
$t_k$ then there is a solution $(x(t),y(t))$ of $f=0$
such that $(x(t_k),y(t_k))=(x_k,y_k)$ for $k=1,\dots,s$.
\enddemo

One can reformulate the theorem to specify not just the value of 
$(x(t),y(t))$ at $t_k$ but also the beginning of its Taylor expansion.  
With  a little more care, the exceptional set $B$ can also be eliminated
(5.1).

\demo{2.5.4 Remark} More generally all of this works if $\c(t)$ is replaced 
with any finite degree extension of
$\c(t)$. These are exactly the fields of meromorphic functions on 
compact Riemann surfaces.
\enddemo

\proclaim{2.6 Low degree equations} $C$ can be described by an equation of
degree at most 2.
\endproclaim

This is of course our starting point, but in 
 higher dimensions this becomes a rather nontrivial question.

It is worthwhile to note the following arithmetic implication: 

\demo{2.6.1 Proposition}  If $\deg f\leq 2$, then $f(x,y)=0$ 
always has a solution over a degree 2 field extension.
\enddemo

In order to see this,   pick $a,b,c$  and consider $f(x,y)=ax+by+c=0$.
Eliminating
$x$ or
$y$  we are left with a quadratic equation in one variable.

\subhead Final remarks about curves
\endsubhead

It should be made clear that the above properties by no means exhaust the
  known characterizations of curves of degree 1 and 2. Some of the others
do not seem to have higher dimensional analogs. I just give a few examples:

\proclaim{ 2.7 Bad characterizations} 
\endproclaim

\demo{2.7.1 Simply connectedness} 

$\pi_1(C(\c))=\{1\}$ iff $\deg f\leq 2$. It turns
out that any smooth hypersurface
$X=(f(x_0,\dots,x_n)=0)\subset
\c\p^n$ is simply connected for $n\geq 3$ \cite{Lefschetz24}. 
\enddemo

\demo{2.7.2 Unique factorization in the coordinate ring} 

The ring $\c[x,y]/f(x,y)$ is a unique factorization domain iff $\deg f\leq 2$.
If $f(x_0,\dots,x_n)$ defines a smooth hypersurface    then 
$\c[x_1,\dots,x_n]/f(x_1,\dots,x_n)$ is a UFD 
for $n\geq 4$ \cite{Grothendieck68}.
\enddemo

\demo{2.7.3 Homogeneous spaces}

If $\deg f=1$ then $C$ is homogeneous under the group $SL(2)$. If $\deg f=2$
then $C$ is homogeneous under the group $O(\bar f)$, the 3-variable
orthogonal group of $\bar f$. In higher dimensions the varieties which are
homogeneous under the action of a linear algebraic group give  
rather special  examples of the class that we want.
\enddemo

\demo{2.7.4 Number of moduli}

Any two lines in $\p^2$ are equivalent under a change of coordinates, and
any two smooth conics in $\c\p^2$ are also equivalent. This fails for $\deg
f\geq 3$.  In all dimensions this property   characterizes hypersurfaces of
degree at most 2, so does not hold for most of the examples in (1.3). (We need
a nondiagonal perturbation to see this.)
\enddemo

The above lists suggest several further possible approaches to low degree
polynomials. Below I list some that do not  to work,  even for curves.

\proclaim{ 2.8 Noncharacterizations} 
\endproclaim

\demo{2.8.1 Topology over $\r$}

One could study curves such that $C(\r)$ is homeomorphic to $S^1$. If $\deg
f\leq 2$ and $C(\r)$ is not empty, this is always the case. Unfortunately,
there are many other curves with this property. For instance, 
 $(x^{2d}+y^{2d}=1)\subset \r\p^2$ is homeomorphic to $S^1$.
\enddemo

\demo{2.8.2 Solutions modulo $p$}

If $f$ has integral coefficients, we can ask about solvability modulo $p$
(or   modulo any number).  

The number of solutions in finite fields are described by the Weil conjectures
(see \cite{Freitag-Kiehl88} for a thorough treatment)   and the degree of $f$
does not affect the asymptotic behaviour much. 
(Though the genus can be computed if we know the exact number of solutions 
modulo $p$ for many values of $p$.)
Low degree equations have
solutions in any finite field \cite{Chevalley35}, but the same holds for
many other cases.
\enddemo

\demo{2.8.3 Solutions in $p$-adic fields}

An equation of degree at most two is not always solvable in $p$-adic fields.
For equations in many variables,
solvability in $p$-adic fields is an interesting question. The rough
picture (which is not quite correct) is that if $f(x_0,\dots,x_n)$ has degree
$d\leq \sqrt{n}$ then
$f$ has a solution in any $p$-adic field and this fails for larger degree.
Thus the
answer does not correspond to our class. 
See \cite{Greenberg69} for a discussion of   these topics.  
\enddemo

\proclaim{ 2.9 Other approaches} 
\endproclaim

\demo{2.9.1 Holomorphic maps $h:\c\to C(\c)$}

If there is a map $\c\p^1\to C(\c)$, then we get plenty of holomorphic
maps $\c\to C(\c)$. If $\deg f\geq 4$ then there are no nonconstant 
holomorphic
maps from $\c$ to $C(\c)$. Unfortunately if $\deg f=3$, then there are
nonconstant holomorphic
maps  $\c\to  C(\c)$.  Thus this property characterizes a slightly
different class of curves. In higher dimensions the two classes differ
substantially. See \cite{Lang86; Vojta91}  for various 
properties of this class.

Vojta pointed out to me that one can consider holomorphic maps $h:\c\to
C(\c)$ whose Nevanlinna characteristic function grows slowly, to get a
characterization of rational curves in the context of the holomorphic theory.
The resulting holomorphic maps are rational, so at the end this is
equivalent to (2.1). 
\enddemo

\demo{2.9.2 The Hasse principle}

One way to overcome the difficulties observed in (1.3) is to refine (1.3.1)
as follows:

Assume that $f(x_1,\dots,x_n)=0$ has a (nontrivial) solution modulo $m$ for
every
$m$ and also over $\r$. Does this imply that $f$ has a solution in $\z$?

(Solvability modulo $m$ for every $m$ is equivalent to solvability in every
$p$-adic field.)

If the answer is yes, one says that the {\it Hasse principle} holds for $f$.
By the Hasse--Minkowski theorem, this is the case 
if $f$ is homogeneous of degree 2.

The  question for higher dimensions is very difficult. It is still not
clear if the  Hasse principle is connected with our class in higher
dimensions or with some smaller class of varieties.
See \cite{Colliot-Th\'el\`ene86,92} for surveys of this direction.
\enddemo

\head 3.  Algebraic Surfaces
\endhead

The next step is to study zero sets of polynomials in three variables
$$
S\deq \{(x,y,z)\vert f(x,y,z)=0\}\subset \a^3.
$$
  It was    noticed in the 19th century that
the true measure of complexity 
of a system of polynomial equations
is the dimension of the set of solutions over $\c$. Thus if
we have 2 equations in 4 variables, the resulting zero set
$$
(f_1(x,y,z,u)=f_2(x,y,z,u)=0)\subset \a^4
$$
behaves to a large extent like surfaces in 3-space. Any surface in 4-space
can be made into a surface in 3-space by a generic projection.
If we generically project a curve in $n$-space to the plane,
the image has only transversal self-intersections. By contrast,
if we  project a surface to 3-space, the image has complicated
self-intersections. According to current view, it is very hard to  study a
surface this way.  (Earlier geometers, being ignorant of this  fact,  proved
rather deep theorems using projections to 3-space.) Thus we are pretty much
forced to look at the general case of varieties:

\demo{Algebraic varieties} Given polynomials $f_1,\dots,f_k$ in $n$
variables, their common zero set
$$
X_{\text{\it aff}}\deq \{(x_1,\dots,x_n)\vert 
f_1(\bold x)=\dots=f_k(\bold x)=0\}\subset \a^n
$$
is called an {\it affine} algebraic variety.  Using homogeneous equations $\bar
f_i$ we obtain {\it projective} varieties
$$
X\deq \{(x_0,\dots,x_n)\vert 
\bar f_1(\bold x)=\dots=\bar f_k(\bold x)=0\}\subset \p^n.
$$
If the coefficients of the $f_i$ are in a field $F$, we  say that $X$ is
{\it defined over} $F$. $X$ is also defined over every bigger field $E\supset
F$, hence $X(E)\subset \text{E}\p^n$, the set of solutions in  $E$,   makes
sense. 
\enddemo

These sets can be very complicated. In order to streamline our discussions, I
make two simplifying assumptions:

{\it All varieties will be irreducible and smooth. }

Over the complex numbers this
means that $X(\c)$ is  a connected manifold. 
These assumptions are satisfied if the coefficients of the $f_i$ are  chosen at
random. The general case can be reduced to this one in various ways. 

The dimension  of  $X$ can be defined in an abstract way. Over $\c$ it  is one
half of the topological dimension of $X(\c)$. This gives the expected value;
for instance if $X\subset \c\p^n$ is defined by a single equation then it has
dimension
$n-1$.

In order to decide which varieties are considered equivalent, we look at the
example of the Mercator projection from (2.1)

\demo{Examples of birational maps} 

(i) Let $S=(x^2+y^2+z^2=1)\subset \r^3$. Project $S$
from the point $(0,0,1)$ to the $(x,y)$-plane $P$.  This provides a one-to-one
map
$$
\pi: S\setminus(0,0,1)\overset{\cong}\to{\longrightarrow} P\cong \r^2.
$$
This looks good, until we notice that projectively there are problems. The
plane is usually compactified as $\r\p^2$, which is not even homeomorphic to
the sphere $S$.

(ii) $H=(x^2-y^2+z^2=1)\subset \r^3$ is a hyperboloid. Project $H$ from the
point
$(0,0,1)$ to the
$(x,y)$-plane $P$.  This provides a one-to-one map
$$
\pi: S\setminus\{(x,y,z)\vert z=1\}
\overset{\cong}\to{\longrightarrow}
 P\setminus\{(x,y)\vert x^2-y^2+1=0\},
$$
and $\pi$ and $\pi^{-1}$ can not be extended to the removed sets in any
reasonable way. Despite this, $\pi$ is clearly very useful in understanding
$S$. For many problems we can use $\pi$ to study 
$S\setminus\{(x,y,z)\vert z=1\}$. The missing set
$\{(x,y,z)\vert z=1\}$ is isomorphic to the plane curve
$\{(x,y)\vert x^2-y^2=0\}$, which is a pair of lines.  

(iii) For $a,b,c\in \q$, $H_{abc}=(ax^2+by^2+cz^2=1)\subset \a^3$ is a
quadric surface. As above, we would like to find a projection of $H$ to a
plane. This can be done over some field, for instance we can project from
$(0,0,1/\sqrt{c})$.  The formulas for $\pi$ and $\pi^{-1}$ involve $\sqrt{c}$,
hence they are of little use if we intend to study $H(\q)$.

If  $a,b,c<0$ then $H_{abc}(\r)$ is empty,  thus there is no map  $g:\r^2\to
H(\r)$. 
\enddemo

\demo{Definition of birational maps}
 Let $X\subset \a^n$ and $Y\subset \a^m$ be affine
varieties. Let $x_i$ (resp. $y_j$) be coordinates on $\a^n$ (resp. $\a^m$).

A {\it rational map}  $g:\a^n\map \a^m$  is given as
$$
g:(x_1,\dots,x_n)\mapsto (g_1(\bold x),\dots,g_m(\bold x)),
$$
where the $g_i$ are rational functions in the variables $x_1,\dots,x_n$. 
Notice that such maps need not be everywhere defined. 
If the coefficients of the $g_i$ are in a field $F$, we say that $g$ is
defined over $F$.

If $g(X)\subset Y$, then  $g$ restricts to a map  $g:X\map Y$.

 We say that
$g:X\map Y$ is {\it birational} if there are subvarieties $A\subsetneq X$ and
$B\subsetneq Y$ such that
 $g$ restricts to a one-to-one map $g:X\setminus A\to Y\setminus B$.

Informally speaking, $X$ and $Y$ are birational if they are isomorphic up to
lower dimensional varieties.

Rational maps of projective varieties can be defined similarly. We can
mimic the above definitions with projective coordinates (in which case the
$g_i$ have to be homogeneous). 
\enddemo
  
A  general introduction to algebraic geometry  can be found in
\cite{Shafarevich94; Hartshorne77}. The analytic theory can be found in
\cite{Wells73; Griffiths-Harris78}. 
The books \cite{Beauville78; BPV84} are
devoted to algebraic surfaces. 
The topological aspects are discussed in
\cite{Donaldson-Kronheimer90; Friedman-Morgan94}.

\subhead  Characterizations of ``low degree" surfaces
\endsubhead

Let $S\subset \p^n$ be a projective surface defined by   homogeneous
equations $f_1=\dots=f_k=0$. For simplicity we always assume
that
$S$ is smooth and   connected.

For surfaces, algebraic geometry provides the basic definition.
Our task is   to see to what extent the other variants (2.2--6) can be
generalized to give an equivalent condition.

\proclaim{3.1  Algebraic geometry} $S$ is rational over $\c$.
\endproclaim

The precise definition of rational is the following:

\demo{3.1.1 Definition} Let $S$ be a smooth projective surface defined over
$\c$.    We say that $S$ is  {\it rational} if there is a
birational map
$g:\c\p^2\map S(\c)$.

If $S$ is defined   over a subfield $F\subset \c$,   we say
that $S$ is {\it rational over $F$} if there is  a
birational map $g:\p^2\map S$ defined over $F$. 
\enddemo

Historically this definition appeared as a rather hard theorem. There are
three classes of surfaces which are  very similar to
rational surfaces, but it is not obvious that they are indeed rational.
These three classes are:

(3.1.2.1) cubic surfaces $S_3\subset \p^3$;

(3.1.2.2) surfaces $S$ which admit a map $f:S\to \p^1$ whose general fiber is
$\p^1$; 

(3.1.2.3)  surfaces  which are images of maps $h:\p^2\map \p^n$.

Cubic surfaces were shown to be rational by \cite{Clebsch1866}. The second
case was settled in \cite{Noether1871} and the third class was treated in
\cite{Castelnuovo1894}.

\proclaim{3.2  Topology}  Homeomorphism versus diffeomorphism.
\endproclaim

Understanding algebraic surfaces in terms of their
topology  turned out to be extremely difficult. 

Some classical questions  can be interpreted in 
topological terms, but this may have been first explicitly done in
\cite{Hirzebruch54}. One of the simplest problems is to give a topological
characterization of the complex projective plane.
This was finally done in  \cite{Yau77}:

\demo{3.2.1 Theorem} Assume that 
$S(\c)$ is homeomorphic to
$\c\p^2$. Then  $S$ is also isomorphic to $\c\p^2$.
\enddemo

The difficulties of this very special case discouraged attempts
to move further in this direction.

A fundamental problem in general is that a birational map $g:S_1\map S_2$
does not induce a homeomorphism.  This question
 can be understood in
terms of the connected sum operation  as follows:

\demo{3.2.2 Proposition} If $S_1(\c)$ and  $S_2(\c)$ are birational then there
are natural numbers
$r,s$ such that
$$
S_1(\c)\#(\overline{\c\p}^2)^r\qtq{is diffeomorphic to} 
S_2(\c)\#(\overline{\c\p}^2)^s,
$$
where $\#$ denotes connected sum and $\overline{\c\p}^2$ is $\c\p^2$ with
reversed orientation. We can assume in addition that $\min\{r,s\}\leq 1$ and
even $\min\{r,s\}= 0$ with a few exceptions. 
\enddemo

In particular we obtain:

\demo{3.2.3 Proposition} If $S$ is rational then
 $S(\c)$ is diffeomorphic to
$$
\c\p^2\#(\overline{\c\p}^2)^r\qtq{or to} 
\c\p^1\times \c\p^1.
$$
\enddemo

(It is not hard to see that 
$(\c\p^1\times \c\p^1)\#\overline{\c\p}^2$ is diffeomorphic to
$\c\p^2\#(\overline{\c\p}^2)^2$, that is  why we have only one
series in (3.2.3).)

By analogy with (2.2) one can ask if the  converse   is also true.
It was noticed some time ago  that the answer is no if we use homeomorphism
instead of diffeomorphism
\cite{Dolgachev66}. As Donaldson theory started to discover the difference
between diffeomorphism and homeomorphism in real dimension 4, the hope emerged
that the converse of (3.2.3) holds for diffeomorphisms. 

This has been one of the motivating questions of the differential topology of
algebraic surfaces. After many contributions, the final step was accomplished
by \cite{Pidstrigach95; Friedman-Qin95}. With the new methods of Seiberg-Witten
theory, the proof is actually quite short \cite{Okonek-Teleman95}:

\demo{3.2.4 Theorem}  Let $S$  be a smooth, projective algebraic surface over
$\c$. Then 
$$
S
\text{ is rational } \qquad\Leftrightarrow \qquad
\matrix
 S(\c) \text{ is diffeomorphic to }\\
\c\p^2\#(\overline{\c\p}^2)^r\text{ or } 
\c\p^1\times \c\p^1.
\endmatrix
$$
\enddemo

\proclaim{3.3  Hard Arithmetic}  $S(\q)$ is ``large".
\endproclaim

Let $S$ be a surface  defined over a number field $F$, most frequently $F=\q$.
As for curves, for any number field $E\supset F$ we   define the counting
function
$$
N_E(S,H)\deq \#\{P\in S(E)\subset \text{E}\p^n\vert 
H(P)\leq H\}.
$$
We hope that $S$ is  rational over $\c$ iff $N_E(S,H)$ grows as a power of
$H$  for some $E$.

Unfortunately this is not quite correct, and there are two related problems. 

(3.3.1.1) Look at the surface $T\deq (x^d+y^d=z^d+u^d)\subset \p^3$.
One can check that $T(\c)$ is smooth. $T$ has high degree, but
$N_{\q}(T,H)$ grows quadratically with $H$.  A closer inspection reveals that
this growth is caused by   
(finitely many) lines on the surface (for instance $(x-z=y-u=0)$)
which contain many rational points. If we remove these lines, there
are very few rational solutions  left.

(3.3.1.2) The growth rate of $N_E(T,H)$ is not
a birational invariant of $T$. Here again the problems are caused by finitely
many rational curves on $T$. 

The examples suggest that we should refine the hope as follows:

\demo{3.3.2 Conjecture}  \cite{FMT89; Batyrev-Manin90} If $T$ is  rational
(over $\c$) then  there is a number field $E$, $0<\beta\in \q$ and
$r\in \n$
such that  
$$
N_E(T\setminus A,H)\qtq{is asymptotic to}  \const\cdot  H^{\beta}(\log H)^r 
$$
for
every sufficiently large subvariety $A\subsetneq T$.
\enddemo

It is furthermore conjectured that  $\beta$ and
$r$  are determined by the geometry of $T$ in a simple way
\cite{Batyrev-Manin90}.  (For higher dimensions these refinements are
problematic, see (4.3).)

A weaker form of (3.3.2) is easy:

\demo{3.3.3 Theorem} If $T$ is  rational (over $\c$) then there is a number
field $E$  and  $\epsilon>0$ such
that 
$N_E(T\setminus A,H)>\const\cdot H^{\epsilon}$ 
for every subvariety $A\subsetneq T$.
\enddemo

The converse of (3.3.2--3) is not quite true. The conceptually correct
formulation will be given in (4.3.2--3). For  surfaces the following form
suffices (cf. \cite{FMT89}). 

\demo{3.3.4 Conjecture}   Assume that  (over $\c$) $T$  is  not
rational  and not birational to $C\times \p^1$ where $C$ is an elliptic curve.
Then  for every  number field $E$ and $0<\epsilon$, there is a subvariety 
$A\subsetneq T$
such that  
$$
N_E(T\setminus A,H) <  \const\cdot  H^{\epsilon}. 
$$
\enddemo

Very little is known in this direction since we have no general
methods to show that nonrational surfaces have only few rational points.

\proclaim{3.4  Complex manifolds} Global holomorphic differential
forms.
\endproclaim

Global holomorphic differential forms on a complex manifold have been much
studied. On a surface we can have $1$-forms and $2$-forms. These are locally
given as 
$$
 f_1dz_1+f_2dz_2, \qtq{respectively}  fdz_1\wedge dz_2,
$$
where  $z_1,z_2$ is a local coordinate system and the $f_i$ are
holomorphic.  In this context, they were first considered in
\cite{Clebsch1868} and systematically studied in
\cite{Picard-Simart1897}.  

As in the curve case, the integrals of these forms over 1- and
2-cycles give basic invariants of a variety \cite{Hodge41}. 
This approach was developed into a very powerful method of studying complex
manifolds, called Hodge theory. 
If there are no global holomorphic differential forms on a surface, then Hodge
theory does not say anything. 

It is easy to see that if $S$ is rational then there are no global
holomorphic differential forms on $S(\c)$. Conversely, one can  hope that
this property characterizes   rational surfaces.

This is 
 close to being
true, and there are two ways of developing a complete answer.

(3.4.1.1) It is known that there are only finitely many families of
exceptions, though the
complete list is not yet known.  

(3.4.1.2) The second approach, which is more promising in higher dimensions,
is to study multivalued
differential forms as well. 
On a surface a multivalued 2-form is locally written as
$f(z_1,z_2)dz_1\wedge dz_2$ where $f$ is a multivalued analytic function.
Thus we may ask about the existence of 2-valued differential forms etc. 
We have the following:

\demo{3.4.2 Theorem}  \cite{Castelnuovo1898} $S$ is rational iff 
there are no global
holomorphic 1-forms  and no global
holomorphic 2-valued 2-forms on $S(\c)$.
\enddemo

It is  technically easier to talk about global sections of
symmetric or tensor powers of the cotangent bundle. In this language the above
result reads:

\demo{3.4.2' Theorem} $S$ is rational iff 
$H^0(S,\Omega_S^1)=0$ and $H^0(S,(\Omega_S^2)^{\otimes 2})=0$.
\enddemo

\proclaim{3.5  Easy Arithmetic}  
There are many solutions over function fields.
\endproclaim

Let $F=\c(t)$  and
$S\subset \text{F}\p^n$ be given by the equations $f_1=\dots=f_k=0$ where
the $f_i$ are homogeneous polynomials in $x_0,\dots,x_n$ with coefficients in
$F$.  Let $\bar F$ denote the algebraic closure of $F$.

The first good news is that the analog of (2.5.2) holds:

\demo{3.5.1 Theorem} \cite{Manin66; Colliot-Th\'el\`ene86}
If $S$ is  rational (over $\bar F$) then $S(F)$ is not
empty.
\enddemo

As for curves, we may want to prove that there are in fact many solutions.
In perfect analogy with (2.5) we have: 

\demo{3.5.2 Theorem}  \cite{KoMiMo92b} Assume that $S$ is  rational (over $\bar
F$).  There is a finite set
$B\subset \c$ such that if $t_1,\dots,t_s\in \c\setminus B$ are arbitrary
points and  $(x_{0k},\dots,x_{nk})$ is any solution of $f_1=\dots=f_k=0$ at
$t_k$ then there is a solution $(x_0(t),\dots,x_n(t))$ of $f=0$
such that $(x_0(t_k),\dots,x_n(t_k))=(x_{0k},\dots,x_{nk})$ for $k=1,\dots,s$.
\enddemo

It would be desirable to generalize to the case when we also
specify  the beginning of the Taylor expansion of $(x_0(t),\dots,x_n(t))$
at certain points. The case when $S$ has a conic bundle structure is quite
easy (see \cite{CTSSD87, I.3.9} for a similar hard arithmetic  proof). The
general case is not known.

All these results hold if $\c(t)$ is replaced  with any
finite degree extension of
$\c(t)$.

\proclaim{3.6 Low degree equations} 
\endproclaim

First we may ask: is every rational surface defined by low degree equations?
The answer is no, there are just too many rational surfaces. It is more
reasonable to ask:

 Is every rational surface $T$ birational to  a surface $S$ which  is
defined by low degree equations?

 By definition, any rational
surface is birational to $\c\p^2$ over $\c$, but
this is rather useless in studying   arithmetic properties of $S$. 
Thus we should be more precise and ask:

\demo{3.6.1 Question}  Let $T$ be a  rational surface defined over a field
$F$.  Is
$T$  always birational over $F$ to  a surface
$S$ which  is defined by low degree equations?
\enddemo

In this form the question is very interesting and fruitful. The answer is
given in two steps.

\demo{3.6.2  Minimal models of surfaces} \cite{Enriques1897}

 The first step is   to simplify the geometry of an arbitrary
smooth projective surface $T(\c)$ by birational maps. The classical name for
this procedure is ``adjunction". Later it was called ``contraction of
(-1)-curves", and the currently fashionable term is ``minimal model program".

 For any surface $T$ we aim to find a birational morphism $f:T\to S$ such that
$S$ is as simple as possible.  (For instance, we may want to make the   Betti
numbers of $S(\c)$   small.) 
$S$ is   called a {\it minimal model} of $T$ (in general it is not unique).

If $T$ is defined over a field $F$, then we can choose $S$ so that
$f$ and $S$ are also defined over  $F$.
(It is  remarkable that the original method of Enriques automatically
works over any field, while the later variants need additional arguments.)
\enddemo

 Next we   study the geometry of the   minimal models $S$ assuming
that $S$ is rational over $\c$. The final result is that there are 4
classes of such surfaces.

\demo{3.6.3 Theorem} \cite{Enriques1897; Manin66; Iskovskikh80c} Let $T$ be a
surface  defined over a field $F\subset \c$ such that $T$ is rational over
$\c$. Then any minimal model of $T$ falls in one of four classes.
(For simplicity, I   use affine coordinates.)

\smallpagebreak
(3.6.3.1) (One low degree equation)

\noindent $S=(f(x,y,z)=0)\subset \a^3$ where  $f$ satisfies one of the
 weighted degree conditions:
$$
\alignat2
\deg(x,y,z)=(1,1,1)&&\qtq{and} \deg f\leq 3&\qtq{(e.g.}x^3+y^3+z^3+1); \\
\deg(x,y,z)=(1,1,2)&&\qtq{and} \deg f\leq 4&\qtq{(e.g.}x^4+y^4+z^2+1);\\
\deg(x,y,z)=(1,2,3)&&\qtq{and} \deg f\leq 6&\qtq{(e.g.}x^6+y^3+z^2+1).\\
\endalignat
$$

\smallpagebreak

(3.6.3.2) (Two low degree equations)

\noindent
$S=(f_1(x,y,z,u)=f_2(x,y,z,u)=0)\subset \a^4$ where $\deg f_i=2$.

\smallpagebreak

(3.6.3.3) (Two equations with low degree in certain variables)

\noindent
$S=(f_1(x,y)=f_2(x,y,z,u)=0)\subset \a^4$ where $\deg f_1=2$ and the degree
of $f_2$ in the $(z,u)$ variables is 2. (The degree of $f_2$ in the $(x,y)$
variables can be high.) 

In these three cases a general choice of  $f,f_1,f_2$ always gives a rational
surface.

\smallpagebreak

(3.6.3.4)  (Miscellaneous) 

These are inconvenient to pin down with equations. They are all birational to
a surface $S=(f(x,y,z)=0)\subset \a^3$ where  $\deg f\leq 9$, but 
$f$ has to be very special. It is much better to notice that   all these
remaining cases are birational to a homogeneous space under a  linear algebraic
group. 
\enddemo

These results imply the following arithmetic assertion:

\demo{3.6.4 Theorem}  Let $S$ be a  surface defined over a field $F\subset
\c$  which is rational over $\c$. 
Then there is a  field extension $E\supset F$ such that
$\deg[E:F]\leq 9$ and 
$S(E)$ is not empty.
\enddemo

\head 4.  Higher Dimensional Varieties
\endhead

After surfaces, the next step is the study of algebraic threefolds. 
The theory of threefolds is much more complicated than the theory of
surfaces, but in the last 20 years a rather satisfactory approach
to threefolds was developed. We know much less about higher dimensions, but
all the conjectures   predict that higher dimensional varieties behave
exactly like threefolds, although the proofs are unknown to us.

Of course it may happen that a few examples will completely change this
picture, but for the moment there is no point in discussing threefolds and
higher dimensional varieties separately.

In the surface case one can always consider only  irreducible and smooth
surfaces.  Starting with dimension three, the smoothness assumption is too
strong, but this is a technical question  which has very little to do with the
essential points of our discussion.

For simplicity,  I mostly consider  smooth varieties. At a few places,
where singularities do cause trouble, I mention this explicitly.

The aspects of higher dimensional algebraic geometry that are discussed here
are treated in the books \cite{CKM88; Koll\'ar96a}. Some other works dealing
with related topics  are
\cite{Ueno75; Koll\'ar et al.92}.  For symplectic topology see
\cite{McDuff-Salamon94,95}.

\subhead  Characterizations of ``low degree" varieties
\endsubhead

Let $X\subset \p^n$ be a smooth projective variety defined by   homogeneous
equations $f_1=\dots=f_k=0$.

As for surfaces,   the algebraic geometry condition gives the basic concept,
but here it takes some work to establish the correct definition.

\proclaim{4.1  Algebraic geometry} $X(\c)$ is rationally connected.
\endproclaim

Already in the surface case it is not easy to show that all low degree
surfaces are rational.
Therefore it did not come as a big surprise that in higher dimensions
rational varieties are too special. A cubic
hypersurface  $X_3^n\subset \c\p^{n+1}$ certainly has low degree. M.
Noether
 knew  that there is a map $p:\c\p^n\map X_3^n$ which is generically 2:1, 
but  nobody was able to prove that $X_3^n$ is rational for $n\geq 3$. 
(And indeed,  $X_3^3$ is not rational \cite{Clemens-Griffiths72}.)
This
leads
 to the following  notion:

\demo{4.1.1 Definition}  $X$ is {\it unirational} (over $\c$) if there is a
rational map
$p:\c\p^n\map X(\c)$ with dense image, where $n=\dim X$.
\enddemo

Very low
degree hypersurfaces in $\c\p^n$ are unirational \cite{Morin40b}. 
Unfortunately,  it seems that the class of unirational varieties is  still
too restrictive.

A new concept was proposed in \cite{KoMiMo92b}. Instead of trying to emulate
global properties of $\c\p^n$, we concentrate on rational curves. $\c\p^n$
has lots of rational  curves (lines, conics and many higher degree ones).
These are images of maps $\c\p^1\to \c\p^n$. The defining property of the new
class should be the existence of lots of maps $\c\p^1\to \c\p^n$. There are
several a priori ways of making this precise. Fortunately, many of these are
equivalent:

\demo{4.1.2 Theorem} \cite{KoMiMo92b} Let $X$ be a smooth  projective  variety
over $\c$. The following are equivalent:

(4.1.2.1) There is an open subset $\emptyset\neq U\subset X(\c)$ such that for
every
$x_1,x_2\in U$ there is a morphism $f:\c\p^1\to X$ satisfying
$x_1,x_2\in f(\c\p^1)$.

(4.1.2.2) For every
$x_1,x_2\in X(\c)$ there is a morphism $f:\c\p^1\to X$ satisfying
$x_1,x_2\in f(\c\p^1)$.

(4.1.2.3) For every
$x_1,\dots,x_n\in X(\c)$ there is a morphism $f:\c\p^1\to X$ satisfying
$x_1,\dots,x_n\in f(\c\p^1)$.

(4.1.2.4)  Let $p_1,\dots,p_n\in\c\p^1$ be distinct points. For each $i$ let
$f_i:D(p_i)\to X(\c)$ be a holomorphic map from a  small disc around $p_i$ to
$X(\c)$.  Let $n_i$ be natural numbers. Then there is a morphism
$f:\c\p^1\to X$ such that the Taylor series of $f_i$ and of $f|D(p_i)$
coincide up to order $n_i$ for every $i$. 

(4.1.2.5) There is a morphism $f:\c\p^1\to X$
such that $f^*T_X$ is ample (see \cite{ibid} for a definition of ample).
\enddemo

\demo{4.1.3 Definition} A smooth projective variety $X$  is
called {\it rationally connected} if it satisfies the equivalent properties
in (4.1.2).
\enddemo

Thus among $n$-dimensional varieties we have 3 classes, with the following
easy containment relations:
$$
\{\text{rational}\}\subset 
\{\text{unirational}\}\subset 
\{\text{rationally connected}\}.
$$
Much effort went into understanding the precise relationship between these
classes. Since 1910, several authors claimed to have produced examples of
rationally connected but nonrational threefolds, but the first correct proofs
appeared only around 1970. By now the situation is quite satisfactory:

\demo{4.1.4 Examples of  rationally connected varieties which are not
rational}

(4.1.4.1)  Dimension three. 

The first examples were quartic 3-folds $X_4\subset
\c\p^4$ 
\cite{Iskovskikh-Manin71} and cubic 3-folds $X_3\subset
\c\p^4$ \cite{Clemens-Griffiths72}. Further development    by
 \cite{Beauville77; Iskovskikh80b; Bardelli84}  gave a quite
complete picture in dimension three.

(4.1.4.2) Conic bundles.

After some very special examples   \cite{Artin-Mumford72}, a general theory
was developed in
\cite{Sarkisov81,82}. This shows that $X_{d,2}\subset \c\p^n\times \c\p^2$ is
not rational for $d\gg 1$. Further examples are in  \cite{Koll\'ar96b}.

(4.1.4.3) Quadric bundles.

Only some special examples are known \cite{CTO89; Peyre93}.

(4.1.4.4) Hypersurfaces

$X_5\subset
\c\p^5$  is considered in  \cite{Pukhlikov87}; the method should  give all
$X_n\subset\c\p^n$. These techniques also give many more  examples 
as in (1.3), see \cite{CPR96}. 
Very general hypersurfaces 
$X_d\subset \c\p^{n+1}$  for $2n/3+2\leq d\leq n+1$ are treated in
\cite{Koll\'ar95}.

(4.1.4.5)  Hypersurface bundles

$X_{c,d}\subset \c\p^m\times \c\p^{n+1}$ where $c\geq 2m$ and
$2n/3+2\leq d\leq n+1$  are considered in  \cite{Koll\'ar96b}.
\enddemo

As this list suggests, most rationally connected varieties are not rational.
Some of the varieties on the above list are unirational, thus rational and
unirational are indeed different notions.  Despite the long list of settled
cases, there are many open problems. I  mention two about hypersurfaces; they
indicate how little is known.

\demo{4.1.5 Some unsolved cases}

(4.1.5.1) Is the general cubic n-fold 
$X_3^n\subset\c\p^{n+1}$ rational for $n\geq 4$? The case of cubic 4-folds has
received a lot of attention. It is known that some special ones are rational
\cite{Morin40a; Tregub93}.  In particular this would show that  rationality
is not deformation invariant.

(4.1.5.2) Is there any rational  (smooth) hypersurface of degree at least 4?
There is very little evidence either way. 
\enddemo

The biggest unsolved question in this picture is the following:

\demo{4.1.6 Conjecture} 
Most rationally connected varieties are not
unirational.
\enddemo

At the moment, there is not a single example known. 
The simplest case to study may be  general quartic threefolds $X_4\subset
\c\p^4$. 

Assume that $X$ is unirational, that is, there is a map $p:\c\p^n\map X$.
The images of linear subspaces show that through a general point of $x\in X$
there are unirational subvarieties of every dimension. Even this weaker
property may fail in general:

\demo{4.1.7 Question} 

Let $X_d\subset \c\p^n$  be a  hypersurface of degree  $d\leq n$ (thus
$X$ is rationally connected). Is it true that for every point $x\in X$ there
is a rational surface $x\in S_x\subset X$?

It is easy to see that this is the case if $\binom{d+1}{2}\leq n$, and
probably also for slightly larger values of $d$. 

I do not see any obvious way to construct rational surfaces when $d$ is close
to $n$. 
\enddemo

Finally I mention another problem concerning rationally connected
varieties. 

\demo{4.1.8 Conjecture}  Let $f:X\to Z$ be a morphism between smooth
projective varieties. Assume that $Z$ and the general fiber $F$ are
rationally connected. Then  $X$ is rationally connected. 
\enddemo

It is easy to see that the special case when $Z=\p^1$ implies the general
one, thus (4.5.1) implies (4.1.8).

\proclaim{4.2  Topology} Diffeomorphism versus symplectomorphism.
\endproclaim

Guided by the results of the surface case, one can look for
three types of theorems in higher dimension:

\demo{4.2.1 Basic Questions}  

(4.2.1.1) Determine all algebraic varieties of a given
topological type.

(4.2.1.2) Relate the topological properties of birationally equivalent
varieties. 

(4.2.1.3) Characterize rationally connected varieties in terms of their
topology. 
\enddemo

As in (3.2), the best example in the first direction is the following
result of \cite{Hirzebruch-Kodaira57; Yau77}

\demo{4.2.2 Theorem} If $X(\c)$ is homeomorphic to
$\c\p^n$ then $X$ is isomorphic to $\c\p^n$.
\enddemo

There are very few such results known, and the proofs  use
rather lucky coincidences.  One may want to have a more
modest aim in mind, and try to show that the topological
structure of $X(\c)$ determines $X$ up to finite ambiguity.
 I noticed the following special case some time ago
 (a proof is given in (5.3)):

\demo{4.2.3 Theorem} Let $M$ be a compact differentiable
manifold with $\dim H_2(M,\q)=1$. Then there are only
finitely many families of algebraic varieties $X$ such that
$X(\c)$ is diffeomorphic to $M$.
\enddemo

For $M$ arbitrary this no longer  holds. This is already shown by the example
of minimal ruled surfaces, but a more convincing   
negative result was
observed by 
\cite{Friedman-Morgan88b}. This shows that diffeomorphism of algebraic 3-folds
is not as strong as for surfaces:

\demo{4.2.4 Example}  Let $S_i$ be 
smooth projective surfaces such that $S_i(\c)$ is simply
connected. Set $X_i\deq S_i\times \c\p^1$.

 For differentiable
manifolds of  real dimension 6, homeomorphism frequently
implies diffeomorphism \cite{Wall66; Sullivan77; Zubr80}. We find that
if $S_i(\c)$ and $S_j(\c)$ are homeomorphic, then
$X_1(\c)$ and $X_2(\c)$ are even diffeomorphic. This gives  several 
unpleasant examples:

(4.2.4.1) Let $S_1$ be a rational surface which is homeomorphic to
a nonrational surface $S_2$ (3.2).  Then $X_1$ is rational,
hence also rationally connected and $X_2$ is not even rationally
connected.

(4.2.4.2) One can construct infinitely many surfaces
$S_i$ such that the $S_i(\c)$ are all homeomorphic, but the
$S_i$ are quite different as algebraic surfaces
\cite{Okonek-V.d.Ven86; Friedman-Morgan88a}.
Thus the manifolds $X_i(\c)$ are all diffeomorphic, but the varieties $X_i$  do
not fit into finitely many families.
\enddemo

\demo{4.2.5 The Topology of Birational Maps} 

Let $X_1$ and $ X_2$ be smooth projective varieties, birational to each other.
In contrast with the surface case, it is not known how the manifolds $X_1(\c)$
and $X_2(\c)$ are related. There are certain surgery type operations, called
blow-ups, that take the role of connected sum with $\overline{\c\p}^2$.
Unfortunately it is not known whether one can go from 
$X_1(\c)$
to $X_2(\c)$ by repeated application of   blow-ups. This is a  hard
problem.

The minimal model program establishes a class of surgery type operations
that can be used  to go from 
$X_1(\c)$
to $X_2(\c)$. At the moment these operations are not well understood from the
topological point of view. 
Furthermore, the intermediate stages involve singular
topological spaces. In dimension three they are all rational homology manifolds
\cite{Koll\'ar91, 2.1.7}, but even this fails in higher dimensions. 
\enddemo

As example (4.2.4) shows, the diffeomorphism type alone does not
characterize rationally connected varieties.  In order to obtain a suitable
analog of (3.2.4), it is necessary to study an additional structure on
$X(\c)$:

\demo{4.2.6 Symplectic manifolds}

A {\it symplectic} manifold is a pair $(M^{2n},\omega)$ where $M$ is a
differentiable  manifold of dimension $2n$ and $\omega$ is a  2-form
$\omega\in
\Gamma(M,\wedge^2T^*)$ which is $d$-closed and nondegenerate.  That is,
$d\omega=0$ and $\omega^n$ is nowhere zero. 

Any smooth projective variety admits a symplectic structure. This can be
constructed as follows. 
  On $\c^{n+1}$
consider the Fubini--Study 2-form
$$
\omega'\deq 
\frac{\sqrt{-1}}{2\pi}
\left[
\frac{\sum dz_i\wedge d\bar z_i}
{\sum |z_i|^2}
-
\frac{(\sum \bar z_i dz_i)\wedge (\sum  z_i d\bar z_i)}
{\left(\sum |z_i|^2\right)^2}
\right].
$$
It is closed, nondegenerate on $\c^{n+1}\setminus\{0\}$ and invariant under
scalar multiplication. Thus $\omega'$ descends to a symplectic 2-form
$\omega$ on $\c\p^n=(\c^{n+1}\setminus\{0\})/\c^*$. 

If $X\subset \c\p^n$ is any smooth variety, then the restriction
$\omega|X$ makes $X(\c)$  into a symplectic manifold. 

The resulting symplectic manifold $(X(\c), \omega|X)$ depends
on the embedding $X\hookrightarrow \c\p^n$, but the dependence is rather
easy to understand:  

We say that two symplectic manifolds $(M,\omega_0)$ and $(M,\omega_1)$ are
{\it symplectic deformation equivalent} if there is a continuous family of 
symplectic manifolds $(M,\omega_t)$ starting with
$(M,\omega_0)$ and ending with $(M,\omega_1)$.

To every smooth projective variety the above construction associates a
symplectic manifold $(X(\c), \omega|X)$ which is unique up to symplectic
deformation equivalence.
\enddemo

This allows us to formulate the proper generalization of (3.2.4):

\demo{4.2.7 Conjecture} Let $X_0$ and $X_1$ be smooth projective
varieties defined over $\c$ such that $(X_0(\c),\omega_0)$ 
is symplectic deformation
equivalent to $(X_1(\c),\omega_1)$. Then
$X_0$ is rationally connected iff $X_1$ is.
\enddemo

The evidence for this conjecture comes from three sources:

The first thing to check is that (4.2.7) holds if there is a continuous
family of algebraic varieties $\{X_t, t\in [0,1]\}$. This case is 
settled:

\demo{4.2.8 Theorem} \cite{KoMiMo92b, 2.4} Let $\{X_t, t\in [0,1]\}$ be
a continuous
family of smooth projective  varieties. 
Then
$X_0$ is rationally connected iff $X_1$ is.
\enddemo

Second, one should try to analyze the examples (4.2.4). 
This was studied in detail by \cite{Ruan94} who showed that the symplectic
structure of $S\times \c\p^1$ can be used to study the differentiable
structure of $S$ in many cases.

The third piece of  evidence is given by the following
closely related result, whose formulation requires  a definition.

\demo{4.2.9 Definition}  A smooth  projective 
variety $X$ over $\c$ is called {\it uniruled}, if it satisfies the following 
equivalent conditions:

(4.2.9.1) There is an open subset $\emptyset\neq U\subset X(\c)$ such that for
every
$x\in U$ there is a morphism $f:\c\p^1\to X$ satisfying
$x\in f(\c\p^1)$.

(4.2.9.2) For every
$x\in X(\c)$ there is a morphism $f:\c\p^1\to X$ satisfying
$x\in f(\c\p^1)$.
\enddemo

The proof of the next result is outlined in  (5.4):

\demo{4.2.10 Theorem} Let $X_0,X_1$ be smooth projective
varieties defined over $\c$ such that $(X_0(\c),\omega_0)$ 
is symplectic deformation
equivalent to $(X_1(\c),\omega_1)$. Then
$X_0$ is uniruled iff $X_1$ is.
\enddemo

(4.2.7) holds if $\dim H_2(X_0,\q)=1$, since  then $X$ is rationally connected
iff it is uniruled \cite{KoMiMo92a}.

It should be noted that if $X_0$ is Fano (4.6.2.1), $X_1$ need not be Fano, as
shown by the examples of rational ruled surfaces.

It would also be interesting to find  some topological properties of rationally
connected varieties. The only general result is the following:

\demo{4.2.11 Theorem} \cite{Campana91b; KoMiMo92b} Let  $X$ be a rationally
connected  variety. Then $X(\c)$ is simply connected. 
\enddemo

\proclaim{4.3  Hard Arithmetic}  $X(\q)$ is ``large".
\endproclaim

As for surfaces, the guiding principle is the following conjecture,
which is a natural generalization of a problem of \cite{Batyrev-Manin90}.

\demo{4.3.1 Conjecture}   If $X$ is  rationally connected
(over
$\c$) then there  are $r\in \n$, $0<\beta\in\q$ and a number field
$F'\supset F$ such that
$$
N_E(X\setminus A,H)\qtq{is asymptotic to}  \text{const}\cdot H^{\beta}(\log
H)^r
$$
for every sufficiently large subvariety $A\subsetneq X$, and for every number
field $E\supset F'$.
\enddemo

The key point  is that  $\beta$ is positive.
Even the following weaker form is completely open: 

\demo{4.3.1' Conjecture}   If $X$ is  rationally connected
 then there  is an $\epsilon>0$ such that  
$$
N_E(X\setminus A,H)>
H^{\epsilon}\qtq{(for $H\gg1$),}
$$
for every   subvariety $A\subsetneq X$, and for every
sufficiently large number field $E$.
\enddemo

There are many special cases where (4.3.1) holds
\cite{FMT89; Batyrev-Manin90;  Batyrev-Tschinkel95}. 
There is a more precise version of the conjecture \cite{Batyrev-Manin90}
asserting that
   the  numbers $\beta, r$ are computable from the
geometry of $T$.   This has been checked in certain cases, but a recent
example of \cite{Batyrev-Tschinkel96} shows that the conjecture for the value
of
$r$ is incorrect.

A  precise computation of the growth  of the number of integral solutions
of the equations
$$
\align
x_1^3+x_2^3+x_3^3&=y_1^3+y_2^3+y_3^3\\
x_1+x_2+x_3&=y_1+y_2+y_3\\
\endalign
$$
is contained in 
 \cite{Vaughan-Wooley95}. This corresponds to (4.3.1) for a
certain singular cubic threefold.  The results confirm (4.3.1), but they
also seem to contradict the more refined conjecture about $r$. Further special
cases are treated in \cite{EMS96}.

The converse of (4.3.1) again fails, but not by much:

\demo{4.3.2 Conjecture}   Assume that   $X$  is  not uniruled (over $\c$).
Then  for every  number field $E$ and $0<\epsilon$, there is a subvariety 
$A\subsetneq X$
such that  
$$
N_E(X\setminus A,H) <  \const\cdot  H^{\epsilon}. 
$$
\enddemo

\demo {4.3.3 The general case}  The problem for a  general variety $X$ can
be reduced to the above two cases as follows. 

Assuming (4.1.8),  there is a
map
$f:X\map Z$ such that
$Z$ is not uniruled and the fibers of $f$ are rationally connected
\cite{KoMoMi92b}. 

Thus we can 
study the points of $X$ in $E$ in two steps. First we have to find the
$E$-points of
$Z$ using (4.3.2). Then for every  $P\in Z(E)$ we study the $E$-points
in the fiber $f^{-1}(P)$, which is rationally connected.
\enddemo

\proclaim{4.4  Complex manifolds} Global holomorphic differential forms.
\endproclaim

As in the surface case, one can study multivalued global holomorphic
differential forms on  $X(\c)$. It
is easy to see that if $X$ is rationally connected, then there are no such
forms:

\demo{4.4.1 Proposition} Let $X$ be a smooth projective variety over $\c$. 
Assume that
$X$ is  rationally connected. Then
$$
H^0\bigl(X,\bigl(\Omega_X^1\bigr)^{\otimes m}\bigr)=0\qtq{for every $m>0$.}
$$
\enddemo

The converse is  conjectured to be true, but it is known only in dimension
three: 

\demo{4.4.2 Theorem} \cite{KoMiMo92b} Let $X$ be a smooth
projective threefold over $\c$. The following are
equivalent:

(4.4.2.1) $X$ is rationally connected;

(4.4.2.2) $H^0\bigl(X,\bigl(\Omega_X^1\bigr)^{\otimes m}\bigr)=0$ for every
$m>0$.
\enddemo

In contrast with (3.5), the current proofs of (4.4.2) require the
vanishing for all values of $m$. It is quite likely that finitely
many of these values are sufficient, but there is no conjecture for
the precise bound.   \cite{KoMiMo92b}
contains further results in this direction.

\proclaim{4.5  Easy Arithmetic}  
There are many solutions over function fields.
\endproclaim

Let $F=\c(t)$  and
$X\subset \text{F}\p^n$ be a subvariety.  Let $\bar F$ denote the algebraic
closure of
$F$.

The higher dimensional analog of  (3.5.1) is
open:

\demo{4.5.1 Conjecture} If $X$ is  rationally connected then $X(F)$ is not
empty.
\enddemo

This is known in many special instances (see, e.g. \cite{Koll\'ar96a, IV.6}),
but these results give very few hints about  the general case.

This of course means that we are also unable to prove that $X$ has many
points in $F$. Surprisingly, one can prove that if $X(F)$ is not empty,
then it is very large. I formulate  the result in the geometric version,
which is more precise.

\demo{4.5.2 Theorem} \cite{KoMiMo92b, 2.13} Let $X$ be a projective variety
over $\c$  and
$f:X\to C$  a morphism onto a smooth curve. Assume that $f$ has a section
$\sigma:C\to X$. Let $c_1,\dots,c_k\in C$ be closed points such that
$f^{-1}(c_i)$ are smooth and  rationally connected. Pick
arbitrary points $p_i\in f^{-1}(c_i)$.

Then $f$ has a section $s=s_{p_1,\dots,p_k}:C\to X$ such that
$s(c_i)=p_i$ for every $i$.
\enddemo

The following more general version is open. In analogy with the number
theoretic terminology (cf. \cite{Mazur92}), it should be called ``weak
approximation for rationally connected varieties over function fields". 

\demo{4.5.3 Conjecture}  Let $X$ be a smooth projective variety
over $\c$  and
$f:X\to C$  a morphism onto a smooth curve whose general
 fibers are   rationally connected. Let $c_1,\dots,c_k\in
C$ be closed points and  $c_i\in D(c_i)\subset C$ small discs around $c_i$. 
Pick local sections 
 $s_i:D(c_i)\to X$ and natural numbers $n_i$.

Then $f$ has a section $s:C\to X$ such that the Taylor series of  
$s|D(c_i)$ agrees with the Taylor series of $s_i$  up to order $n_i$, for every
$i$.
\enddemo

In the special case when $X=C\times Y$, this follows from (4.1.2.4).

\proclaim{4.6 Low degree equations} 
\endproclaim

As in the surface case, the principal question is the following:

\demo{4.6.1 Question}  Let $X$ be a  rationally connected variety defined over
a field
$F$.  Is
$X$  always birational over $F$ to  a variety
$Y$ which  is defined by low degree equations?
\enddemo

In contrast with the surface case, this is interesting even  for $F=\c$.

In analogy with (3.6), first we need:

\demo{4.6.2.  Minimal model program}

  This is a general method to simplify the
structure of an arbitrary smooth  projective variety. Already in dimension 3 it
is rather complicated
(cf. \cite{Mori82,88}), and in higher dimensions remains conjectural.
See \cite{Koll\'ar87,90} for introductions.
The program can
be performed over any field $F$ with minor modifications.

For rationally connected varieties we end up with a  variety $Y$
(birational to $X$) satisfying one of the following conditions:

(4.6.2.1)  $Y$ is a Fano variety, that is, $-K_Y$ is ample.
Unfortunately, $Y$ may be singular. The singularities are rather
mild (terminal and $\q$-factorial), but they do cause certain problems.

(4.6.2.2)  There is a morphism $p:Y\to Z$ such that $Z$ and the fibers of $p$
are rationally connected.
\enddemo

In the second case we hope  to reduce
problems about $X$ to questions about $Z$ and about the fibers of $f$. Thus we
mainly concentrate on the first case. 
Some of
the basic questions are settled:

\demo{4.6.3 Theorem} (4.6.3.1) \cite{Nadel91; Campana91a; KoMiMo92a,c} For any
$n$ there are only finitely many families of smooth Fano varieties of dimension
$n$.

(4.6.3.2) \cite{Kawamata92} There are only finitely many families of singular
Fano threefolds arising in (4.6.2.1).
\enddemo

In both cases the proof yields explicit (though huge) bounds on the number of
families and also on the degrees of the defining equations of the Fano
varieties.

In dimension three there is a complete list of all smooth Fano varieties, but
no such list exists in the singular case. In any case, classifying Fano
threefolds up to isomorphism may not be the sensible thing to do.
Our original  variety $X$ is determined by $Y$ only up to birational
equivalence; thus it makes sense   to classify rationally connected
threefolds  up to birational equivalence. \cite{Alexev94; Corti96} contain
significant steps in this direction. 

\demo{4.6.4 Listing by low degree equations} 

 Smooth Fano threefolds  were studied by
G. Fano in a series of articles  spanning four decades starting in 1908.
A modern account of these  works was given in  \cite{Iskovskikh80a,b}. 
The results of
\cite{Mukai89} give a  better description, especially over nonclosed fields. 
For singular Fano threefolds there is no general theory; a series of
examples can be found in
\cite{Fletcher89}. 

If there is a morphism $p:X\to Z$ as in (4.6.2.2), then the results of
(3.6) give us
   defining equations as in (3.6.3).
Instead of listing all cases, I just give two examples:

(4.6.4.1) $S=(f_1(u,v)=f_2(x,y,z,u,v)=0)\subset \a^5$, 

\noindent where $\deg f_1=2$ and the degree
of $f_2$ in the $(x,y,z)$ variables satsifies one of the conditions of
(3.6.3.1) (The degree of
$f_2$ in the
$(u,v)$ variables can be high.)

(4.6.4.2)
$S=(f_1(x_1,x_2)=f_2(x_1,\dots,x_4)=f_3(x_1,\dots,x_6)=0)\subset
\a^6$, 

\noindent where $\deg f_1=2$, the degree of $f_2$ in the $(x_3,x_4)$
variables is 2 and the degree of $f_3$ in the $(x_5,x_6)$ variables is 2.
 (The degrees  in the other
 variables can be high.) 

In both cases a general choice of the $f_i$ gives  a rationally connected
variety.
\enddemo

These results imply the following arithmetic consequence:

\demo{4.6.5 Theorem}    There is a constant $D(3)$ with the following property:

Let $X$ be a  rationally connected threefold defined over a field $F\subset
\c$. 
Then there is a  field extension $E\supset F$ such that
$\deg[E:F]\leq D(3)$ and 
$X(E)$ is not empty.
\enddemo

One can write down an explicit bound for $D(3)$, though I have not done it.
Conjecturally, a similar result holds in any dimension.

\head 5. Appendix
\endhead

The aim of this appendix is to outline the proofs of some statements  which
are new or for which I could not find a suitable reference. 

\demo{5.1 Proposition} Let $B$ be a smooth proper curve over $\c$ and
$f:S\to B$ a proper ruled surface. Let $b_i\in B$ be different points
and $D(b_i)$ a small disc around $b_i$.
Let $s_i:D(b_i)\to S$ be holomorphic (or formal) sections and $n_i$ natural
numbers.

Then there is a section $s:B\to S$ such that $s|D(b_i)$ agrees with $s_i$ up
to order $n_i$ for every $i$.
\enddemo

\demop $S$ is birationally trivial; that is, there is a birational map
$\pi: \p^1\times B\map S$. We obtain local sections 
$$
s'_i\deq \pi^{-1}\circ s_i:D(b_i)\to \p^1\times B.
$$
Assume that it takes $k$ blow-ups to resolve the indeterminacies of $\pi$.
Let $s':B\to \p^1\times B$ be a section such that $s'|D(b_i)$ agrees with
$s'_i$ up to order $n_i+k$ for every $i$. Then we can take $s\deq \pi\circ
s'$.

Thus it is sufficient to find $s'$. Equivalently, we need to find a map
$\bar s: B\to \p^1$ with prescribed local behavior $\bar s_i:D(b_i)\to \p^1$.
By a generic coordinate   change in $\p^1$ we can assume that 
$\bar s_i(b_i)\in \c$ for every $i$. 

Choose another point $b_0$. One can always find regular functions on
the affine curve 
$B\setminus\{b_0\}$ with prescribed local behaviour at the points $b_i$.
\qed\enddemo

\demo{5.2 Proof of (4.1.2.4)} 
We need to show that (4.1.2.4) is implied by (4.1.2.3). 
As a first step, I prove the following weaker
version:

(5.2.1)  Let $p_1,\dots,p_n\in\c\p^1$ be disctinct points. For each $i$ let
$f_i:D(p_i)\to X(\c)$ be a holomorphic map from a  small disc around $p_i$ to
$X(\c)$.  Let $n_i$ be natural numbers. Then there is a morphism
$g:\c\p^1\to X$  and holomorphic maps $h_i:D(p_i)\to  \c\p^1$ such that the
Taylor series of $f_i$ and of
$g\circ h_i|D(p_i)$ coincide up to order $n_i$ for every $i$.

To see this, let $D\subset \c$ be the unit disc and
$f,g:D\to \c^n$ two holomorphic maps with coordinate functions $f^j, g^j$.
Assume that $f(0)=g(0)=0\in \c^n$.
Let $B_0\c^n\to \c^n$ be the blow-up of $0\in \c^n$. $f$ and $g$ lift to
holomorphic maps $\bar f,\bar g:D\to B_0\c^n$. Explicit local computation
shows the following: 

(5.2.2.1) If $\bar f$ and $\bar g$ agree up to order $n$, then so do $f$ and
$g$. 

(5.2.2.2) If $f^1(t)=g^1(t)=t$ and 
$\bar f$ and $\bar g$ agree up to order $n-1$, then
$f$ and $g$ agree up to order $n$.

Using (5.2.2.1)  for repeated blow-ups, we  first reduce (5.2.1) to the case
when the $f_i$ are immersions. 
Then up to  a local coordinate change
 we may assume that $f_i^1(t)=t$ for every $i$. 
We can now prove (5.2.1) by induction on $\sum n_i$, since (4.1.2.3) gives it
for $\sum n_i=0$.

The only subtle point is the reduction step from order $1$ to order $0$.
Let $p\in D\subset \c\p^1$ be a disc. 
Given an immersion $f:D\to X$,  let $x=f(p)$ and  $\pi:B_xX\to X$ be the
blow-up with exceptional divisor $E\subset B_xX$.  Assume that we have $\bar
g:\c\p^1\to B_xX$ such that 
$\bar f$ and $\bar g$ agree up to order $0$ at $p$. We would like to conclude
that $f$ and $g\deq \bar g\circ \pi$ agree up to order $1$ at $p$. (5.2.2.2)
gives this, if
$g$ is an immersion. Thus we have to choose $\bar g:\c\p^1\to B_xX$
to be transversal to $E$. This is slightly stronger than (4.1.2.3), but can
easily be arranged (see the proofs of II.3.14 and IV.3.9 in
\cite{Koll\'ar96a}). 

Once we have (5.2.1), we just need to find a  map $h:\c\p^1\to \c\p^1$
which approximates every $h_i$ up to order $n_i$ and set $f\deq g\circ h$

The $f$ we found is a multiple cover of a curve in $X$. 
As in \cite{Koll\'ar96a, IV.3.9} we can perturb  $f$ to obtain another solution
of (4.1.2.4) where $f|\c\p^1\setminus\{p_1,\dots,p_n\}$ is an embedding.
\qed\enddemo

 \demo{5.3 Proof of (4.2.3)} 
Assume that $X(\c)$ is diffeomorphic to $M$.
 We use 
the  formula  \cite{Hirzebruch66,20.3.6*}
$$
\chi(\o_X) = \sum_{s\geq 0}  
\frac1{2^{n+2s}(n-2s)!}c_1(X)^{n-2s}A_s(p_1,\dots,p_s)[M],\tag 5.3.1
$$
where the $A_s$ are certain polynomials of the Pontrjagin classes of $M$
and $A_0=1$.
From Hodge theory we know that
$$
|\chi(\o_X)|\leq \sum \dim_{\c} H^i(X,\o_X)\leq \sum \dim_{\c} H^i(M,\c),
$$
and so $\chi(\o_X)$ is bounded in terms of $M$. Since 
$b_2(M) =1 $,
we can fix an ample divisor $H$ in $\pic(X)$ and then 
$c_1(X)\equiv rH$ for some rational number $r$. (5.3.1) becomes a polynomial
 equation for $r$. As $\chi(\o_X)$ runs through all the
possible values, we get only finitely many possible values for $r$. 
 Therefore the
self-intersection number $(H^n)$ and the intersection number $(c_1(X)\cdot
H^{n-1})$ are
bounded depending on $M$ only. The result now follows from  Matsusaka's Big
Theorem (in the form given in \cite{Koll\'ar-Matsusaka83}).\qed\enddemo

The proof provides an  effective bound on the number of 
families of algebraic structures on a given manifold $M$,
but the bound is
enormous even in the simplest cases.

\demo{5.4 Proof of (4.2.10)} The proof is an application of the theory of
Gromov--Witten invariants. I recall the main concepts in the needed special
case. See
\cite{McDuff-Salamon94,95} for details of the general theory.

Let $X$ be a smooth projective variety over $\c$. Fix a point $x\in X$,  
a homology class $A\in
H_2(X(\c),\z)$ and  very
ample divisors in general position $H_i\subset X$, $i=1,\dots, k$.  

Let $y_0,\dots,y_k\in \c\p^1$ be general points. For suitable   $k$,
there may be only finitely many maps
$$
f:\c\p^1\to X\qtq{such that} 
f_*[\c\p^1]=A, f(y_0)=x,\text{ and }  f(y_i)\in H_i,\ i=1,\dots, k.
$$
We define an invariant
$$
\tilde F_{A,X}(x,H_1,\dots,H_k; y_0,\dots,y_k)\deq \text{ the number of such
maps.}
\tag 5.4.1 
$$

Gromov's theory of pseudo-holomorphic curves shows that one can make a
similar definition where $X$ is replaced by a symplectic manifold
$(M,\omega)$ endowed with a  general almost complex structure.
The corresponding invariant is denoted by
$$
\tilde \Phi_{A,M,\omega}(x,H_1,\dots,H_k; y_0,\dots,y_k).\tag 5.4.2
$$
It is one of the  Gromov--Witten invariants of $(M,\omega)$. 
In fact, this is an invariant of the symplectic deformation equivalence class.

In general the algebraic number (5.4.1) and the symplectic number (5.4.2) are
different.  Under suitable conditions they are equal, and this means that we
can get information about rational curves on $X$ from the symplectic structure
$(X(\c),\omega_X)$ (4.2.6).  
This idea was used by \cite{Ruan93} to show that the  extremal rays of
Mori theory can be described using the symplectic structure.
We need the following two results.
(In \cite{Ruan93} they are proved under the extra assumption that the
symplectic structure is semi-positive. This is no longer necessary.)

\demo{5.4.3 Theorem} 
Let $X$ be a smooth projective variety over $\c$ and $(M,\omega)$ the
corresponding symplectic manifold.

(5.4.3.1) If $\tilde \Phi_{A,M,\omega}(x,H_1,\dots,H_k;
y_0,\dots,y_k)\neq 0$, then there is a rational map
$f:\c\p^1\to X$ such that $f_*[\c\p^1]=A$, $f(y_0)=x$ and $f(y_i)\in H_i,\
i=1,\dots, k$. 

(5.4.3.2) 
$\tilde F_{A,X}(x,H_1,\dots,H_k; y_0,\dots,y_k)=
\tilde \Phi_{A,X(\c),\omega_X}(x,H_1,\dots,H_k; y_0,\dots,y_k)$ 
if the
following conditions are satisfied:

(5.4.3.2.1)  If  $g:\c\p^1\to X$ is any map such that 
$g_*[\c\p^1]=A$ and 
$g(y_0)=x$, then $H^1(\c\p^1, g^*T_X)=0$.

(5.4.3.2.2) If    $C_1,\dots,C_m\subset X$ are rational curves such that
$\sum [C_i]=A$ and $x\in C_1$, then $m=1$.
\enddemo

We can now prove (4.2.10). 

Let $(M,\omega)$ be the common symplectic structure of $X_0$ and of $X_1$.
Fix a very general point $x\in X_0$. 
Fix a very ample divisor $H\subset X_0$ and 
let $x\in C\subset X$ be a rational curve
such that $(C\cdot H)$ is minimal ($C$ exists since $X_0$ is uniruled). Set
$A\deq [C]$.  By \cite{KoMiMo92c,1.1}, the
condition (5.4.3.2.1) holds and (5.4.3.2.2) follows from the minimality  of
$(C\cdot H)$. Let $k$ be the dimension of the space of maps
 $g:\c\p^1\to X$  such that 
$g_*[\c\p^1]=A$ and 
$g(y_0)=x$. Let $H_1,\dots, H_k\subset X_0$ be general divisors linearly
equivalent to $H$. By construction,
$\tilde F_{A,X}(x,H_1,\dots,H_k; y_0,\dots,y_k)$ is defined and is nonzero.
Thus  $\tilde \Phi_{A,M,\omega }(x,H_1,\dots,H_k; y_0,\dots,y_k)\neq 0$.

By (5.4.3.1) this implies that there is a rational curve through any
very general point of $X_1$, and thus $X_1$ is also uniruled.\qed\enddemo

Finally we prove that condition (1.2)  correctly identifies the class of
rationally connected varieties among diagonal hypersurfaces.

\demo{5.5 Proposition}  
Let $X$ be any smooth compactification of the affine
hypersurface 
$$
(\sum_{i=1}^n c_ix_i^{d_i}+c_0=0)\subset \c^n.
$$

(5.5.1) $X$ is rationally connected iff $\sum 1/d_i\geq 1$.

(5.5.2) The Kodaira dimension of $X$ is nonnegative iff $\sum 1/d_i< 1$.
\enddemo

\demop  Consider first the case $n=2$, assuming $d_1\leq d_2$. View $X$ as a
$d_1$-sheeted cover of the line ramified along $c_2x_2^{d_2}+c_0=0$. The
Hurwitz formula gives that
$$
2g(X)=(d_1-1)(d_2-2)+(\text{ramification at infinity}).
$$
This implies (5.5) for $n=2$.

If $n\geq 3$ then as in (1.3), we
view these as hypersurfaces in weighted projective spaces.   
Let $d=lcm(d_i)$, $d=d_ia_i$ and set $a_0=1$, $d_0=d$.  A (nonsmooth)
compactification is given by 
 the projective weighted hypersurface
$$
Y:=\sum_{i=0}^n c_ix_i^{d_i}\subset \p(a_0,\dots,a_n).
$$
As long as $\prod c_i\neq 0$, these hypersurfaces are isomorphic (over $\c$),
thus $Y$ can be viewed as a general member of the linear system
$|x_0^{d_0},\dots, x_n^{d_n}|$. This implies that $Y$ has only quotient
singularities and Picard number 1 for $n\geq 4$.

Assume that $d< \sum a_i$. $K_Y=\o(d-\sum a_i)$, thus $Y$ is 
$\q$-Fano. Therefore $Y$ is uniruled by \cite{Miyaoka-Mori86}. Let
$p:\bar Y\to Y$ be a desingularization and  $\bar f:\bar Y^0\to Z$  the MRC
fibration
\cite{KoMiMo92b}. The fibers of $p$ are all rationally connected (cf.
\cite{Koll\'ar96a, VI.1.6.2}), thus $\bar f$ descends to  $f:Y^0\to Z$.
If $n\geq 4$, then as in \cite{Koll\'ar96a,IV.4.14}, we obtain that $Z$ is a
point, hence
$Y$ is rationally connected. If $n=3$ then we use that
$h^1(X,\o_X)=h^1(Y,\o_Y)=0$. A smooth uniruled surface $S$ with 
$h^1(S,\o_S)=0$ is rational, hence $X$ is rational.

Next assume that $d\geq  \sum a_i$. Let $e=d-\sum a_i$ and introduce
$e$ new coordinates  $x_{n+1},\dots,x_{n+e}$ of weight 
$a_{n+1}=\dots=a_{n+e}=1$. Consider the
hypersurface
$$
Z:=\sum_{i=0}^n c_ix_i^{d_i}+\sum_{j=1}^e c_{n+j}x_{n+j}^{d}
\subset \p(a_0,\dots,a_{n+e}).
$$
Here every $a_i$ divides $d=\sum a_i$, thus $\p(a_0,\dots,a_{n+e})$ has only
index one canonical singularities. Therefore the same holds for $Z$. But
$\omega_Z\cong \o_Z$, and this implies that $\kappa(Z)=0$.

General  fibers of the projection map
$$
Z\map \p^e\qtq{given by} (x_0,\dots,x_{n+e})\mapsto
(x_0,x_{n+1},\dots,x_{n+e})
$$
are isomorphic to  $Y$. This shows that $\kappa(Y)\geq 0$.

Since a variety  can not be rationally connected and have nonnegative Kodaira
dimension at the same time, this proves (5.5).\qed\enddemo

\demo{5.5.3 Remark}  It is not true that $X$ is of general type if
$d>\sum a_i$. For instance,
$x_1^2+x_2^3+\sum c_ix_i^{d_i}+c_0=0$ has an elliptic fiberspace structure
(projection to the $(x_3,\dots,x_n)$-subspace) for every value of the $d_i$.
\enddemo

\demo{Acknowledgements} 
I would like to thank S. Gersten, K. Ribet, Y. Ruan  and P. Vojta for useful
conversations and   e-mails. J.-L. Colliot-Th\'el\`ene sent very detailed
comments which helped me to understand the arithmetical questions much better.
 Partial financial support was provided by the
NSF under grant number  DMS-9102866.
These notes were typeset by \AmSTeX, the \TeX\ macro system of the American
Mathematical Society.
\enddemo

\vskip1cm

University of Utah, Salt Lake City UT 84112 

kollar\@{}math.utah.edu

\end